%% file: main.tex
\documentclass{easychair}

\usepackage{lineno}
\usepackage{doc}
\usepackage{booktabs}
\usepackage{amsmath,amsfonts}
\usepackage{algorithm}
\usepackage{algpseudocode}
\usepackage{graphicx}
\usepackage{multirow}
\usepackage{multicol}
\usepackage{amssymb}
\usepackage{subcaption}
\usepackage{xcolor}
\usepackage{diagbox}
\usepackage{enumitem}
\setlist[itemize]{leftmargin=*}
\usepackage{url}
\usepackage{threeparttable}
\usepackage{pifont}
\usepackage{nopageno}
\usepackage{bm}
\usepackage{scalerel}
\usepackage{makecell}
\usepackage[colorinlistoftodos]{todonotes}

\newcommand{\R}{\mathbb{R}}

\title{
A low-rank balanced truncation approach for large-scale RLCk model order reduction based on extended Krylov subspace and a frequency-aware convergence criterion
\\
}

\author{
Christos Giamouzis,
Dimitrios Garyfallou,
Nestor Evmorfopoulos,\\
and George Stamoulis
}

\institute{
 Dept. of Electrical and Computer Engineering, 
University of Thessaly, Volos, Greece \\
  \email{\{cgiamouzis, digaryfa, nestevmo, georges\}@e-ce.uth.gr}
}

\begin{document}
\maketitle

\input{abstract}
\input{1_Introduction}
\input{2_Background}
\input{3_Proposed}

\input{4_Experimental_Evaluation}

\input{5_Conclusions}

\bibliographystyle{IEEEtran}
\bibliography{easychair}

\end{document}

%% file: abstract.tex
\begin{abstract}
Model order reduction (MOR) is essential in integrated circuit design, particularly~when dealing with large-scale electromagnetic models extracted from complex designs. The numerous passive elements introduced in these models pose significant challenges in the simulation process. MOR methods based on balanced truncation (BT) help address these challenges by producing compact reduced-order models (ROMs) that preserve the original model's input-output port behavior. In this work, we present an extended Krylov subspace-based BT approach with a frequency-aware convergence criterion and efficient~implementation techniques for reducing large-scale models. Experimental results indicate that our method generates accurate and compact ROMs while achieving up to $\times$22 smaller ROMs with similar accuracy compared to ANSYS~RaptorX\texttrademark\ ROMs for large-scale benchmarks.
\end{abstract}

%% file: 1_Introduction.tex
\section{Introduction}
\label{sec:intro}

Electromagnetic model extraction is crucial for designing and verifying integrated circuits (ICs), enabling precise simulation of the passive elements of the design. However, simulating extracted RLCk models with millions of elements and multiple ports is extremely computationally expensive. Model order reduction (MOR) can reduce the complexity of such models while maintaining accurate input/output port behavior~\cite{Odabasioglu1998, Antoulas2004}. By constructing reduced-order models (ROMs) that capture the essential dynamics of the original system, MOR can significantly reduce simulation time, enabling faster design iterations in IC development.

There are two main approaches to MOR. Moment matching (MM) methods are preferred for their efficiency,  but they require manual selection of the number of moments~\cite{Odabasioglu1998}. Most~importantly, they correlate the final ROM size with the number of moments and ports, limiting scalability. On the contrary, balanced truncation (BT) provides explicit theoretical bounds~for~the approximation error and is independent of the number of ports~\cite{Antoulas2004}. However, BT is restricted to small-scale models due to the high computational complexity of solving Lyapunov equations~\cite{Antoulas2004}.

In this paper, we introduce an efficient low-rank BT technique to address the main scalability issue of the conventional BT approach. Specifically, we employ the extended Krylov subspace (EKS) method, which effectively solves the Lyapunov equations, drastically reducing the computational load of BT~\cite{Arxiv_23}. Additionally, we incorporate a frequency-aware convergence criterion, ensuring accuracy in the frequency range of interest. Experimental evaluation indicates that the proposed method can be integrated into commercial extraction tools, such as the ANSYS RaptorX\texttrademark ~\cite{raptX}, to generate more compact ROMs of large-scale multi-port~RLCk~models. 

%% file: 2_Background.tex
\section{Background}
\label{sec:background}

Consider the modified nodal analysis (MNA) description~\cite{MNA} of an $n$-node, $m$-branch (inductive), $p$-input, and $q$-output RLCk circuit in the time domain: 
\begin{equation} \label{mna}
    \begin{aligned}
        \begin{pmatrix}
            \mathbf{G_{n}} & \mathbf{E} \\
            \mathbf{-E}^T & \mathbf{0} 
        \end{pmatrix}
        \begin{pmatrix}
            \mathbf{v}(t) \\
            \mathbf{i}(t) 
        \end{pmatrix} +
        \begin{pmatrix}
            \mathbf{C_{n}} & \mathbf{0} \\
            \mathbf{0} & \mathbf{M} 
        \end{pmatrix}
        \begin{pmatrix}
            \dot{\mathbf{v}}(t) \\
            \dot{\mathbf{i}}(t) 
        \end{pmatrix} = 
        \begin{pmatrix}
            \mathbf{B}_1 \\
            \mathbf{0} 
        \end{pmatrix}
        \mathbf{u}(t), \quad 
        \mathbf{y}(t) = 
        \begin{pmatrix}
            \mathbf{L}_1 \quad
            \mathbf{0} 
        \end{pmatrix}
        \begin{pmatrix}
            \mathbf{v}(t) \\
            \mathbf{i}(t) 
        \end{pmatrix},
    \end{aligned}
\end{equation}
where $\mathbf{G_{n}} \in\R^{n\times n}$ (node conductance matrix),  $\mathbf{C_{n}}\in\R^{n\times n}$ (node capacitance matrix), $\mathbf{M}\in\R^{m\times m}$ (branch inductance matrix), $\mathbf{E}\in\R^{n \times m}$ (node-to-branch incidence matrix), $\mathbf{v}\in\R^{n}$ (vector of node voltages), $\mathbf{i}\in\R^{m}$ (vector of inductive branch currents), $\mathbf{u} \in \R^{p} $ (vector of input excitations), $\mathbf{B}_1\in\R^{n\times p}$ (input-to-node connectivity matrix), $\mathbf{y}\in\R^{q}$ (vector of output measurements), and $\mathbf{L}_1\in\R^{q\times n}$ (node-to-output connectivity matrix). 
Moreover, we denote
$\dot{\mathbf{v}}(t) \equiv \frac{d\mathbf{v}(t)}{dt}$ and $\dot{\mathbf{i}}(t) \equiv \frac{d\mathbf{i}(t)}{dt}$.
If we now define the model order as $N \equiv n + m$, the state vector as $\mathbf{x}(t) \equiv \begin{pmatrix}
\mathbf{v}(t) \\
\mathbf{i}(t) 
\end{pmatrix}$, and~also:
\begin{equation*}
\begin{aligned}
\mathbf{G}\equiv -\begin{pmatrix}
\mathbf{G_{n}} & \mathbf{E} \\
\mathbf{-E}^T & \mathbf{0} 
\end{pmatrix},\quad  \mathbf{C} \equiv \begin{pmatrix}
\mathbf{C_{n}} & \mathbf{0} \\
\mathbf{0} & \mathbf{M} 
\end{pmatrix},\quad \mathbf{B}\equiv \begin{pmatrix}
\mathbf{B}_1 \\
\mathbf{0} 
\end{pmatrix}, \quad \mathbf{L}\equiv \begin{pmatrix}
\mathbf{L}_1 \quad
\mathbf{0} 
\end{pmatrix}
\end{aligned}, 
\end{equation*}
then Eq.~(\ref{mna}) can be written in the generalized state-space form, or so-called descriptor form:
\begin{equation}
\begin{aligned} \label{state}
\mathbf{C}\frac{d \mathbf{x}(t)}{d t} = \mathbf{G x}(t) + \mathbf{Bu}(t), \quad 
\mathbf{y}(t) = \mathbf{L x}(t).
\end{aligned}
\end{equation}

\noindent The objective of MOR is to produce an equivalent ROM: 
\begin{equation}
\begin{aligned}
\mathbf{\tilde C} \frac{d \mathbf{\tilde x}(t)}{d t} =\mathbf{\tilde G} \mathbf{\tilde x}(t) + \mathbf{\tilde B} \mathbf{u(t)}, \quad
\mathbf{\tilde y}(t) = \mathbf{\tilde L \tilde x}(t),
\end{aligned}
\end{equation}				
where  $\mathbf{\tilde G}, \mathbf{\tilde C} \in \R^{r\times r} $, $\mathbf{\tilde B} \in \R^{r\times p} $, $\mathbf{\tilde L} \in \R^{q\times r}$, the reduced order $r<<N$, and the output error is bounded as $||\mathbf{\tilde{y} }(t) -\mathbf{y}(t)||_2 < \varepsilon||\mathbf{u}(t)||_2$ for given  $\mathbf{u}(t)$ and small $\varepsilon$. The output error bound can be expressed in the frequency domain as $||\mathbf{\tilde{y} }(s) -\mathbf{y}(s)||_2 < \varepsilon||\mathbf{u}(s)||_2$ via Plancherel's theorem~\cite{Plancherel}.~If
\begin{equation*}
\begin{aligned}
\mathbf{H}(s) = \mathbf{L}(s\mathbf{C} - \mathbf{G})^{-1} \mathbf{B}, \quad
\mathbf{\tilde H}(s) =  \mathbf{\tilde L}(s\mathbf{\tilde C} - \mathbf{\tilde G})^{-1} \mathbf{\tilde B}
\end{aligned}
\end{equation*}
are the transfer functions of the original model and the ROM, the corresponding
output
error~is: 
\begin{equation*}
\begin{aligned}
||\mathbf{\tilde{y} }(s) -\mathbf{y}(s)||_2 = ||\mathbf{\tilde{H}}(s) \mathbf{u}(s) - \mathbf{H}(s)\mathbf{u}(s)||_2 \quad \leq \quad ||\mathbf{\tilde{H}}(s) - \mathbf{H}(s)||_\infty||\mathbf{u}(s)||_2,
\end{aligned}
\end{equation*}
where $||.||_\infty$ is the 
$\mathcal{L}_2$ matrix norm or $\mathcal{H}_\infty$ norm of a rational transfer function. 
Thus, to~bound this error, we need to bound the distance between the transfer functions:~$||\mathbf{\tilde{H}}(s) - \mathbf{H}(s)||_\infty~<~\varepsilon$.

%% file: 3_Proposed.tex
\section{MOR by Balanced Truncation}
\label{sec:bt}

BT relies on the computation of the controllability Gramian $\mathbf{P}$ and observability Gramian~$\mathbf{Q}$, which are calculated as the solutions of the following Lyapunov matrix equations \cite{Antoulas2004}:
\begin{equation}
\begin{aligned}\label{Eq:lyap_sta}
(\mathbf{C}^{-1}\mathbf{G}) \mathbf{P} +  \mathbf{P} (\mathbf{C}^{-1}\mathbf{G})^T  = - (\mathbf{C}^{-1}\mathbf{B}) (\mathbf{C^{-1}}\mathbf{B})^T, \quad
(\mathbf{C}^{-1}\mathbf{G})^T \mathbf{Q} +  \mathbf{Q}(\mathbf{C}^{-1}\mathbf{G}) = - \mathbf{L}^T \mathbf{L}.
\end{aligned}
\end{equation}

The controllability Gramian $\mathbf{P}$ describes the degree to which the states are controllable by the inputs, while the observability Gramian $\mathbf{Q}$ reflects the degree to which the states are observable at the outputs. A ROM can theoretically be generated by eliminating the states that are difficult to control or observe. However, in the original state-space coordinates, certain states may be easy to control but difficult to observe, and vice versa. The process of ``balancing'' transforms the state vector to a new coordinate system, where the controllability and observability of each state are balanced, meaning each state is equally difficult to control and observe. An appropriate transformation $\mathbf{Tx}(t)$ exists, leading to the balanced state-space model:
\begin{equation*}
\begin{aligned}
\mathbf{TCT}^{-1}\frac{d (\mathbf{Tx}(t))}{d t} = \mathbf{TGT^{-1}} (\mathbf{Tx}(t)) + \mathbf{TBu}(t), \quad
\mathbf{y}(t) = \mathbf{LT^{-1}} (\mathbf{Tx}(t)).
\end{aligned}
\end{equation*}
This balanced representation preserves the system's transfer function $\mathbf{H}(s)$ and simplifies to $\mathbf{P}$\ =\ $\mathbf{Q}$\ =\ $diag(\sigma_1, \sigma_2, \dots, \sigma_N)$~\cite{Antoulas2004},
~where $\sigma_i$ are the Hankel singular values (HSVs). These HSVs are the square roots of the eigenvalues of the product $\mathbf{P Q}$,  
i.e., $\sigma_i = \sqrt{\lambda_i(\mathbf{P Q})}$.
In the above balanced model,
the states with the largest HSVs are the easiest to both control and observe. If $r$ of them are retained (truncating the $N-r$ states associated with the smallest HSVs), the error between the original and the reduced-order transfer functions is bounded as:
\begin{equation*}
||\mathbf{H}(s) - \mathbf{\tilde H}(s)||_\infty \leq 2( \sigma_{r+1} +\sigma_{r+2}+...+\sigma_{N}).
\end{equation*}
The above serves as an ``a-priori'' criterion that offers flexibility by allowing either the specification of a ROM size $r$ to compute the error or a target error ($target\_error$) to determine~the~number $r$ of HSVs to be preserved. This adaptability is a key advantage of BT over MM methods.

\begin{figure}[!hbt]
\vspace{-1.5em}
\begin{algorithm}[H]
\footnotesize
	\caption{MOR by balanced truncation}\label{vanilla_BT}
    \textbf{Inputs:}  $ \mathbf{G}, \mathbf{C}, \mathbf{B}, \mathbf{L} $\\
	\textbf{Outputs:} $\mathbf{\tilde G}, \mathbf{\tilde C}, \mathbf{\tilde B}, \mathbf{\tilde L}$
	\begin{algorithmic}[1]
		\State Solve the Lyapunov equations to obtain the Gramian matrices $\mathbf{P}$ and $\mathbf{Q}$~\cite{Lathauwer2004}
		\State Compute the
            SVD of the Gramian matrices: $\mathbf{P} = \mathbf{U}_P \mathbf{\Sigma}_P \mathbf{V}_P^T $ and $\mathbf{Q} = \mathbf{U}_Q \mathbf{\Sigma}_Q \mathbf{V}_Q^T $
            \State Find the square root of the Gramian matrices: $\mathbf{Z}_P = \mathbf{U}_P \mathbf{\Sigma}_P^{1/2}$ and $\mathbf{Z}_Q = \mathbf{U}_Q \mathbf{\Sigma}_Q^{1/2}$
            \State Compute the SVD of the product of the roots: $\mathbf{Z}_Q^T\mathbf{Z}_P = \mathbf{U}\mathbf{\Sigma}\mathbf{V}^T$	
            \State Compute transformation matrices: $\mathbf{T}_{(r\times N)}$\ =\ $\mathbf{\Sigma}_{(r\times r)}^{-1/2} \mathbf{U}_{(r\times N)} \mathbf{Z}_Q^T $,\ $\mathbf{T}_{(N\times r)}^{-1}$\ =\ $ \mathbf{Z}_P\mathbf{V}_{(N\times r)}\mathbf{\Sigma}^{-1/2}_{(r\times r)}$
            \State Compute ROM: $\mathbf{\tilde G}$\ =\ $\mathbf{T}_{(r\times N)}\mathbf{G}\mathbf{T}_{(N\times r)}^{-1}$, \ $\mathbf{\tilde C}$\ =\ $\mathbf{T}_{(r\times N)}\mathbf{C}\mathbf{T}_{(N\times r)}^{-1}$,\  $\mathbf{\tilde B}$\ =\ $\mathbf{T}_{(r\times N)}\mathbf{B},\ \mathbf{\tilde L}$\ =\ $\mathbf{L}\mathbf{T}_{(N\times r)}^{-1}$
	\end{algorithmic}
\end{algorithm}
\vspace{-1.5em}
\end{figure}

The main steps of the BT procedure are summarized in Algorithm \ref{vanilla_BT}. The main limitation of BT is its high computational and memory cost, which makes it impractical for large-scale models (with $N$ over a few thousand states).
This is due to the computationally expensive operations required, such as solving Lyapunov equations and performing singular value decomposition (SVD), both of which have a complexity of $O(N^3)$. Additionally, they are applied on dense matrices, since the Gramians $\mathbf{P}, \mathbf{Q}$ are dense even if the system matrices $\mathbf{C}, \mathbf{G}, \mathbf{B}, \mathbf{L}$~are~sparse. 

However, the products $(\mathbf{C}^{-1}\mathbf{B}) (\mathbf{C^{-1}}\mathbf{B})^T$ and $\mathbf{L}^T\mathbf{L}$ have a much lower numerical rank compared~to~$N$, as $p,q<<N$. This results in low-rank Gramians that can be approximated using low-rank techniques, significantly reducing the complexity and memory requirements for solving the Lyapunov equations and performing SVD, which are now performed with a complexity of order $k$ rather than $N$.

\subsection{Low-rank BT MOR}
\label{sec:bt_eks}
The essence of low-rank BT MOR is to iteratively project the Lyapunov equations onto a lower-dimensional Krylov subspace and solve the resulting small-scale equations to obtain low-rank approximate solutions of Eq.~(\ref{Eq:lyap_sta}). The $k$-dimensional standard Krylov subspace is defined as:

\begin{equation*} 
    \mathcal{K}_k(\mathbf{G}_{C},\mathbf{B}_{C}) = span \{\mathbf{B}_{C},\mathbf{G}_{C}\mathbf{B}_{C},  \mathbf{G}_{C}^{2}\mathbf{B}_{C},\dots,\mathbf{G}_{C}^{k-1}\mathbf{B}_{C}\},
\end{equation*}
where $\mathbf{G}_{C} \equiv \mathbf{C}^{-1}\mathbf{G},\hspace{0.25em} \mathbf{B}_{C} \equiv \mathbf{C}^{-1}\mathbf{B}$. 
If $\mathbf{K} \in \R^{N \times k}$ ($k<<N$) is a projection matrix whose columns span the $k$-dimensional standard Krylov subspace, then the projected Lyapunov equation (for the controllability Gramian $\mathbf{P}$) onto $\mathcal{K}_k(\mathbf{G}_{C},\mathbf{B}_{C})$ is:
\begin{equation} \label{Eq:projected_lyapunov_equations}
    (\mathbf{K}^T\mathbf{G}_{C}\mathbf{K})\mathbf{X} +\mathbf{X} (\mathbf{K}^T\mathbf{G}_{C}\mathbf{K})^T =- \mathbf{K}^T\mathbf{B}_{C}\mathbf{B}_{C}^T\mathbf{K}
\end{equation}
(the same holds true for the observability Gramian $\mathbf{Q}$ with $\mathbf{G}_{C}^T$, $\mathbf{L}^{T}$ in place of $\mathbf{G}_{C}$, $\mathbf{B}_{C}$). The solution $\mathbf{X} \in \R^{k \times k}$ of  Eq.~(\ref{Eq:projected_lyapunov_equations}) can be back-projected to the $N$-dimensional space to give an approximate solution $\mathbf{P} \approx \mathbf{K}\mathbf{XK}^T$ for the original large-scale Eq.~(\ref{Eq:lyap_sta}), and a low-rank factor $\mathbf{Z} \in \R^{N \times k}$ of $\mathbf{P}$ can be obtained as $\mathbf{Z} =\mathbf{K}\mathbf{U}\mathbf{\Sigma}^{1/2}$, where $[\mathbf{U},\mathbf{\Sigma},\mathbf{V}] = SVD(\mathbf{X})$ and $\mathbf{P} \approx \mathbf{Z}\mathbf{Z}^T$.

While the projection process is independent of the chosen subspace, its effectiveness heavily relies on it.
The convergence to an accurate solution can be accelerated by enhancing the standard Krylov subspace $\mathcal{K}_k(\mathbf{G}_{C},\mathbf{B}_{C}) $ with information from the subspace $\mathcal{K}_k(\mathbf{G}_{C}^{-1},\mathbf{B}_{C})$, which corresponds to the inverse matrix $\mathbf{G}_{C}^{-1}$, leading to the EKS~\cite{Arxiv_23, ASPDAC21}:
\begin{equation} \label{Eq:eksm}
\mathcal{K}_k^C(\mathbf{G}_{C},\mathbf{B}_{C}) = 
span \{\mathbf{B}_{C},  \mathbf{G}_{C}^{-1}\mathbf{B}_{C}, \mathbf{G}_{C}\mathbf{B}_{C},\mathbf{G}_{C}^{-2}\mathbf{B}_{C}, \mathbf{G}_{C}^{2}\mathbf{B}_{C},\dots, \\
\mathbf{G}_{C}^{-(k-1)}\mathbf{B}_{C}, \mathbf{G}_{C}^{k-1}\mathbf{B}_{C}\}.
\end{equation}
The EKS method (EKSM) begins with the vectors $\{ \mathbf{B}_{C}, \mathbf{G}_{C}^{-1}\mathbf{B}_{C} \}$ and iteratively builds an EKS $\mathcal{K}_k^C(\mathbf{G}_{C},\mathbf{B}_{C})$ of increasing dimension, solving the projected Lyapunov Eq.~(\ref{Eq:projected_lyapunov_equations}) in each iteration, until a sufficiently accurate approximation of the solution of Eq.~(\ref{Eq:lyap_sta}) is achieved. The complete EKSM is presented in Algorithm~\ref{eksm-alg}. 
Below are some efficient implementation details:
\begin{itemize}
\itemsep0em
    \item \textbf{Matrix inversion by linear solves}: 
	Algorithm \ref{eksm-alg} uses the system matrices $\mathbf{G}$, $\mathbf{C}$ or $\mathbf{G}^T$, $\mathbf{C}^T$ instead of $\mathbf{G}_{C} \equiv \mathbf{C}^{-1}\mathbf{G}$ or $\mathbf{G}^T_{C} \equiv (\mathbf{C}^{-1}\mathbf{G})^T$ since the (generally dense) inverse matrices are only required for products with $p$ vectors (in step ~\hyperref[eksm-alg:step2]{2}) and $2pj$ vectors (in steps ~\hyperref[eksm-alg:step4]{4} and ~\hyperref[eksm-alg:step11]{11} of each iteration), which can be handled as linear solves like $\mathbf{C}\mathbf{Y} = \mathbf{R}$ and $\mathbf{G}\mathbf{Y} = \mathbf{R}$ (or $\mathbf{C}^T\mathbf{Y} = \mathbf{R}$, $\mathbf{G}^T\mathbf{Y} = \mathbf{R}$), using either direct or iterative methods~\cite{Bavier2012}.
  
  \item \textbf{Handling of sparse/dense matrices}: Matrix $\mathbf{M}$ of Eq.~(\ref{mna}) is typically very dense due to the huge number of mutual inductances. To efficiently handle both sparse ($\mathbf{C}_n$) and dense ($\mathbf{M}$) matrix blocks of $\mathbf{C}$, we use 
  specialized data structures and numerical~techniques. This includes parallel CPU-optimized methods for sparse matrices and GPU-accelerated techniques~\cite{CMG_GPU} for dense matrices.
  \item \textbf{Solution of the small-scale Lyapunov equations}: 
  To solve the small-scale~($2pj \times 2pj$) Lyapunov~equations  in step~\hyperref[eksm-alg:step5]{5} of each iteration,
  we employ the Bartels-Stewart~algorithm~\cite{Lathauwer2004}. 
  \item  \textbf{Convergence criterion:} The error estimation~\cite{SMACD2024} relies on the ROM transfer function $\mathbf{\tilde H}(s)$ and is described by:
\begin{equation*}
    \max_{i=1 \dots l}{\frac{||\mathbf{\tilde H}_{j}(s_i)-\mathbf{\tilde H}_{j-1}(s_i)||_\infty}{||\mathbf{\tilde H}_{j}(s_i)||_\infty}},
\end{equation*}
where $\mathbf{\tilde H}_{j}(s_i)$ is the ROM transfer function at the $j$-th iteration (calculated at frequency $s_i = 2\pi f_i$) and $l$ is the number of evaluated frequency points evenly distributed across a frequency range [$f_{min}$, $f_{max}$]. The proposed criterion offers insight into the extent to which the transfer function changes between iterations at the frequencies of interest. Moreover, it proves to be practical and effective for circuit simulation problems, where designers are only interested in the circuit's behavior in certain frequency windows. The iterative procedure stops when the error remains below a certain threshold ($tol$) for three consecutive iterations.
\end{itemize}
\begin{figure}[!hbt]
\vspace{-1.5em}
\begin{algorithm}[H]
\footnotesize
\caption{Extended Krylov subspace method for low-rank solution of Lyapunov equations}\label{eksm-alg}	
	\textbf{Input:}  $ \mathbf{G}_{C} \equiv \mathbf{C}^{-1}\mathbf{G},  \mathbf{B}_{C} \equiv \mathbf{C}^{-1}\mathbf{B}$  (or  $\mathbf{G}_{C}^T$, $\mathbf{L}^T $)\\
	\textbf{Output:} $\mathbf{Z}$ such that $\mathbf{P} \approx \mathbf{Z} \mathbf{Z} ^T $
	\begin{algorithmic} [1]
		\State $j=1$; $p=size\_col(\mathbf{B}_{C})$
        \label{eksm-alg:step2}
		\State $\mathbf{K}^{(j)} = Orth([\mathbf{B}_{C},\mathbf{G}_{C}^{-1}\mathbf{B}_{C}])$
		\While {$j < maxiter$}
        \label{eksm-alg:step4}
    		\State $\mathbf{A} = \mathbf{K}^{(j)T}\mathbf{G}_{C}\mathbf{K}^{(j)} $;\quad $\mathbf{R} = \mathbf{K}^{(j)T}\mathbf{B}_{C} $  
        \label{eksm-alg:step5}
    		\State Solve $\mathbf{A}\mathbf{X} +\mathbf{X} \mathbf{A}^{T}  = -\mathbf{R}\mathbf{R}^{T}$ for  $\mathbf{X} \in \R^{2pj \times 2pj}$
    		\If {converged}
        		\State $[\mathbf{U}, \mathbf{\Sigma}, \mathbf{V}] = \mathbf{SVD}(\mathbf{X})$; \quad
        		$\mathbf{Z} = \mathbf{K}^{(j)}\mathbf{U}\mathbf{\Sigma}^{1/2}$
        	    \State \textbf{break}
        	\EndIf
    		\State $k_1=2p(j-1)$; $k_2 = k_1+p$; $k_3 = 2pj$ 
        \label{eksm-alg:step11}
    		\State $\mathbf{K}_1  = [\mathbf{G}_{C}\mathbf{K}^{(j)}(:,k_1+1:k_2),\mathbf{G}_{C}^{-1}\mathbf{K}^{(j)}(:,k_2+1:k_3)]$
    		\State $\mathbf{K}_2 = Orth(\mathbf{K}_1) $ \quad w.r.t. \quad $\mathbf{K}^{(j)}$
    		\State $\mathbf{K}_3 = Orth(\mathbf{K}_2) $ 
    		\State $\mathbf{K}^{(j+1)} = [\mathbf{K}^{(j)},\mathbf{K}_3]$\quad 
		    \State $j=j+1$
		\EndWhile
	\end{algorithmic}
\end{algorithm}
\vspace{-1.5em}
\end{figure}

%% file: 4_Experimental_Evaluation.tex
\vspace{-1em}
\section{Experimental Evaluation}
\label{sec:experimental}

We evaluated EKSM using RLCk models extracted from various circuits via ANSYS RaptorX\texttrademark~\cite{raptX}. The evaluated designs consist of a phase-locked loop (PLL), an analog mixer, a time-interleaved digital-to-analog converter (TI\_DAC), an injection-locked frequency multiplier (ILFM), a VGA circuit, hybrid couplers (HCs), Wilkinson power dividers (WPDs), and typical transceiver blocks, such as low-noise-amplifiers (LNAs) and oscillators (VCO). Their detailed characteristics are listed in Tables~\ref{table:experimental} and~\ref{table:low_rank_results}. Two experiments were conducted: in the first one, we used small-scale benchmarks ($<$ 30K nodes), where the original and ROM transfer functions could be directly compared; in the second one, we used large-scale benchmarks and compared EKSM to golden RaptorX\texttrademark\,ROMs through S-parameter plotting. For the reduction process, $target\_error$ and $tol$ were set to 1e-2 and the number of frequencies $l$ was set to 20. Experiments were performed on a Linux server with a 2.80 GHz 16-thread CPU and 64 GB of memory.

\begin{table}[hbt!]
\centering
    \caption{Detailed characteristics of RLCk models and experimental results of EKSM}
        \vspace{-0.4em}	
\label{table:experimental}
\footnotesize
\setlength{\arrayrulewidth}{0.2mm}
\setlength{\tabcolsep}{2.75pt}
\renewcommand{\arraystretch}{1.2}
\begin{tabular}{|c||c|c|c|c|c||c|c|c|c|}
\hline
 \multirow{2}{*}{Model} & \multirow{2}{*}{Order} & \multirow{2}{*}{\#\,ports} & \multirow{2}{*}{\makecell{\#\,mutual \\ induct.}} & \multirow{2}{*}{$|| \mathbf{G}_C ||_F$} & \multirow{2}{*}{$cond(\mathbf{G}_C)$} & \multirow{2}{*}{\makecell{ROM \\ order}} & \multirow{2}{*}{Error} & \multirow{2}{*}{\makecell{Reduction \\ time}} & \multirow{2}{*}{Memory} \\ 
 & & & & & & & & & \\ \hline
 PLL\,@\,28\,GHz & 1474 & 4 & 251680 & 7.03e+21 & 4.41e+16 & 96 & 1.21e-03 & 3.8 s & 72 MB \\ \hline
 Mixer\,@\,28\,GHz & 1498 & 10 & 79794 & 4.13e+23 & 1.81e+18 & 100 & 2.04e-04 & 2.09 s & 58 MB \\ \hline
 TI\_DAC\,@\,28\,GHz & 3869 & 160 & 365494 & 4.16e+22 & 9.62e+15 & 1280 & 1.14e-06 & 3 min & 1.18 GB \\ \hline
 LNA\,@\,56\,GHz & 4274 & 6 & 1988882 & 1.73e+24 & 1.37e+19 & 204 & 7.98e-04 & 31 s & 341 MB \\ \hline
 LNA\,@\,28\,GHz & 6956 & 6 & 5360490 & 1.55e+23 & 3.15e+18 & 144 & 4.48e-03 & 30 s & 594 MB \\ \hline
 ILFM\,@\,14\,GHz & 15665 & 11 & 18394794 & 4.68e+27 & 9.02e+22 & 176 & 1.19e-03 & 1 min & 1.54 GB \\ \hline
 LNA\,@\,2.4\,GHz & 25602 & 6 & 72959220 & 6.23e+23 & 5.14e+19 & 144 & 1.42e-03 &  3 min & 6.62 GB \\ \hline
\end{tabular}
\vspace{-1.2em}
\end{table}

For the first experiment, the results are presented in Table~\ref{table:experimental}, where the error refers to the max relative error between the transfer functions of the original models and ROMs, which is calculated as $||\mathbf{H}(s) - \mathbf{\tilde H}(s)||_\infty \ / \ ||\mathbf{H}(s)||_\infty$ at the designated frequencies. For every benchmark, a base frequency of 100 MHz is chosen ($f_{min}$ = 1e+8) and the maximum frequency is set to twice the resonance frequency of each circuit (e.g., $f_{max}$ = 56e+9 for PLL\,@\,28\,GHz). As can be seen, EKSM generates accurate and compact ROMs across every type of benchmark with a maximum error below 0.14\%. Additionally, the convergence criterion effectively strikes a balance between error and final ROM size while being computationally efficient. This is also visible through the performance results, where the reduction time remains below 3 min for benchmarks with less than 30K nodes and the memory requirements are not significantly high.

\begin{figure}[!tbh]
  \centering
    \vspace{-2.55em}
    \includegraphics[width=0.70\textwidth]{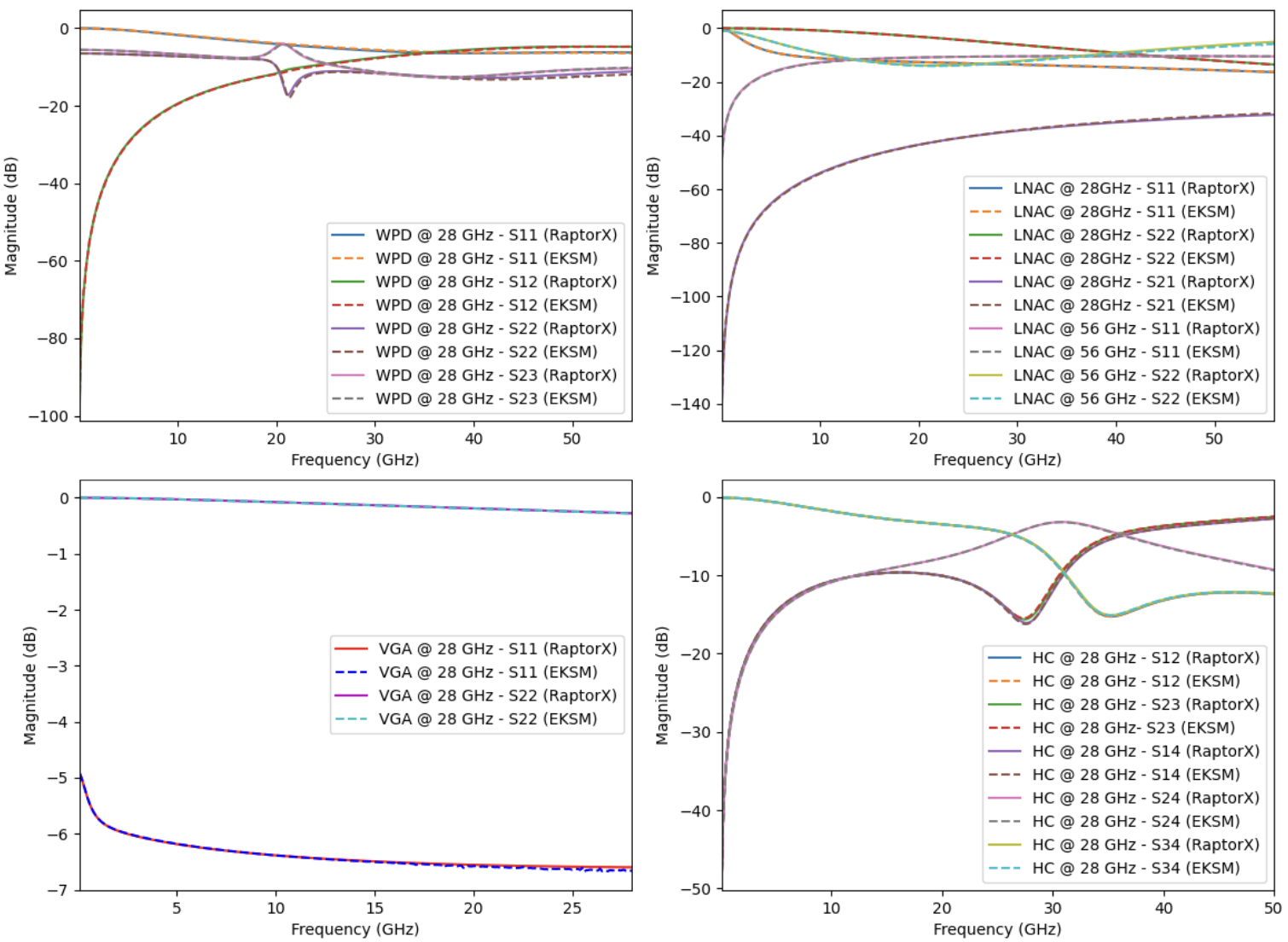}
    \vspace{-8pt}
    \caption{Comparison of accuracy between EKSM and RaptorX\texttrademark~ROMs.}
    \label{fig:accuracy_results}
    \vspace{-1.5em}
\end{figure}

For the second experiment, the results are demonstrated in Table~\ref{table:low_rank_results}. The S-parameters plots of Figure~\ref{fig:accuracy_results} indicate that EKSM achieves accuracy close to that of RaptorX\texttrademark\ while producing on average $\times$13.2 more compact ROMs. Although EKSM has higher reduction~time and memory requirements, they are still reasonable and can be significantly improved in futurk work.

\begin{table}[!hbt]
\centering
    \caption{ ROM order and MOR performance of EKSM vs RaptorX\texttrademark}
              \vspace{-0.4em}	

    \label{table:low_rank_results}
\footnotesize
\setlength{\arrayrulewidth}{0.2mm}
\setlength{\tabcolsep}{2.75pt}
\renewcommand{\arraystretch}{1.2}
\begin{tabular}{|c||c|c|c||c|c||c|c||c|c|}
\hline
\multirow{2}{*}{Model} & \multirow{2}{*}{\makecell{Initial \\ order}} & \multirow{2}{*}{\#ports} & \multirow{2}{*}{\makecell{\#mutual \\ induct.}} & \multicolumn{2}{c||}{ROM order} & \multicolumn{2}{c||}{Reduction time} & \multicolumn{2}{c|}{Memory (GB)} \\ 
\cline{5-10}& & & & RaptorX\texttrademark& EKSM & RaptorX\texttrademark& EKSM & RaptorX\texttrademark& EKSM\\ \hline
VGA\,@\,28\,GHz & 95189 & 13 & 126766838 & 4744 & 286 & 1 min & 11 min & 32.63 & 19.14 \\ \hline
HC\,@\,56\,GHz & 98024 & 5 & 165802476 & 1267 & 120 & 2 min & 8 min & 24.05 & 25.05 \\ \hline
WPD\,@\,56\,GHz & 100888 & 4 & 193641938 & 765 & 40 & 3 min & 4 min & 24.79 & 29.76 \\ \hline
VCO\,@\,13\,GHz & 104367 & 4 & 188436057 & 407 & 56 & 2 min & 4 min & 26.48 & 28.96 \\ \hline
LNAC\,@\,56\,GHz & 128574 & 9 & 169339965 & 2172 & 378 & 1 min & 19 min & 25.82 & 26.01 \\ \hline
WPD\,@\,28\,GHz & 129087 & 4 & 259462454 & 885 & 40 & 3 min & 5 min & 25.35 & 34.57 \\ \hline
HC\,@\,28\,GHz & 134710 & 5 & 264162513 & 787 & 90 & 4 min & 8 min & 24.31 & 35.62 \\ \hline
LNAC\,@\,28\,GHz & 162881 & 11 & 323090671 & 4768 & 308 & 6 min & 27 min & 78.52 & 49.58 \\ \hline
\end{tabular}
\vspace{-0.6em}
\end{table}

%% file: 5_Conclusions.tex
\section{Conclusions}
\label{sec:conclusions}
An alternative MOR technique for accurately reducing large-scale RLCk models is introduced. The proposed low-rank BT approach incorporates an iterative EKS projection method with a frequency-aware convergence criterion to produce accurate and compact ROMs. Experimental results demonstrate that our method provides up to ×22 smaller ROMs than ROMs obtained by~ANSYS~RaptorX\texttrademark\,\,for large-scale benchmarks with negligible deviations in S-parameters.